\documentclass[a4paper, 10pt]{article}
\usepackage[utf8]{inputenc}
\usepackage[T1]{fontenc}

\usepackage{amsmath}
\usepackage{amsfonts}
\usepackage{amsthm}
\usepackage{microtype}
\usepackage{hyperref}
\usepackage[pass]{geometry}
\usepackage{cite}

\usepackage{mathrsfs}
\usepackage{fouriernc}

\usepackage{multirow}
\usepackage{array}
\usepackage{booktabs}
\usepackage{pgfplotstable}
\usepackage{colortbl}
\usepackage{pdflscape}

\usepackage{rotating}

\pgfplotstableset{
    math/.style={preproc cell content/.append style
                ={/pgfplots/table/@cell content/.add={$}{$}}},
    col sep = semicolon,
    every even row/.style={before row={\rowcolor[gray]{0.95}}},
    every head row/.style={before row=\toprule,after row=\midrule},
    every last row/.style={after row=\bottomrule},
    row/5/.append ={after row=\bottomrule},
    columns/Group/.style={
        string type,
        column type=l, 
        column name={},
        math
    },
    columns/S/.style={
        int detect,
        column type=r,
        column name={$|S|$}
    },
    columns/Size/.style={
        int detect,
        column type=r,
        column name={Size}
    },
    columns/d/.style={
        int detect,
        column type=r,
        column name={$d$}
    },
    columns/B2d/.style={
        int detect,
        column type=r,
        column name={$|B_{2d}|$}
    },
    columns/n/.style={
        int detect,
        column type=r,
        column name={$n$}
    }, 
    columns/m/.style={
        int detect,
        column type=r,
        column name={$m$}
    },
    columns/t-pre/.style={
        sci,
        sci sep align,
        precision=1,
        column name={$t_\textsl{pre}\;\;[s]$}
    },
    columns/lambda0/.style={
        dec sep align,
        precision=5,
        column name={$\lambda_0$}
    },
    columns/(1-d)lambda0/.style={
        dec sep align,
        fixed,
        fixed zerofill,
        precision=5,
        column name={$(1-\delta)\lambda_0$}
    },
    columns/t-SDP/.style={sci,
        sci sep align,
        precision=1,
        column name={$t_\textsl{SDP}\;\;[s]$}
    },
    columns/prec/.style={sci,
        sci sep align,
        precision=1,
        column name={\textsl{prec}}
    },
    columns/resid/.style={
        sci,
        sci sep align,
        sci zerofill,
        precision=1,
        column name={$\|r\|_1<$}
    },
    columns/kappa-a/.style={
        fixed,
        precision=5,
        dec sep align,
        column name={$\kappa_a$}
    },
    columns/kappa-n/.style={
        fixed,
        fixed zerofill,
        precision=4,
        dec sep align,
        column name={$\kappa_n$}
    },
    columns/kappa-our/.style={
        fixed,
        fixed zerofill,
        precision=5,
        dec sep align,
        column type/.add={|}{|},column name={$\kappa$}
    }
}

\newcommand{\Hilbert}{\mathscr{H}}
\newcommand{\SL}{\operatorname{SL}}
\newcommand{\supp}{\operatorname{supp}}
\newcommand{\Z}{\mathbb{Z}}

\newtheorem{theorem}{Theorem}
\newtheorem{lemma}[theorem]{Lemma}

\title{Certifying numerical estimates\\of spectral gaps}
\author{Marek Kaluba\thanks{Author has been parially supported by the National Science Centre, under grant 2015/19/B/ST1/01458}~~and Piotr  W. Nowak\thanks{This project has received funding from the European Research Council (ERC) under the European Union's Horizon 2020 research and innovation programme (grant agreement no. 677120-INDEX)}}

\begin{document}
\maketitle

\begin{abstract}
We establish a lower bound on the spectral gap of the Laplace operator on special linear groups using conic optimisation. In particular, this provides a constructive (but computer assisted) proof that these groups have Kazhdan property (T). 
Software for such optimisation for other finitely presented groups is provided.
\end{abstract}

\section{Introduction}
Estimation of spectral gaps and Kazhdan constants for finitely presented groups is a fundamental problem in analytic group theory. 
In particular, much work was devoted only to establishing and then improving bounds for $\SL(n,\Z )$, see \cite{Burger1991, Shalom1999, Kassabov2005, Hadad2007}. 
These estimates are analytical in nature, often providing a homogeneous bound for all $n$. In this paper we obtain significantly 
better numerical bounds for low dimensional special linear groups ($n=3,4,5$) as well as for special linear groups over certain finite fields.

It is known that property (T) is equivalent to the positivity of the operator $\Delta^2 - \lambda \Delta$ in the maximal 
group $C^*$-algebra, for some positive $\lambda$. This, in turn, is equivalent, by the Positivstellensatz \cite{Schmudgen2009}, to an approximation of $\Delta^2 - \lambda \Delta$ by sums-of-(hermitian)-squares, i.e. elements of the form $\sum \xi_i^*\xi_i$ \cite{Netzer2013}.
Ozawa in \cite{Ozawa2016} showed that this can be effectively checked, as the decomposition is attainable in $\mathbb{Q}[G]$, hence in a finite-dimensional subspace of the real group ring $\mathbb{R}[G]$. Moreover, he proved that the approximation error term can be modified to sum-of-squares, at the expense of $\lambda$.

In \cite{Netzer2015} Netzer and Thom cast the problem into the language of semidefinite (numerical) optimisation and, as a proof-of-concept, numerically improved the existing bound for the spectral gap of $\Delta$ (and thus for the Kazhdan constant) for $\SL(3,\Z )$ by three orders of magnitude.
To achieve this they provide a constructive estimation of residuals which inevitably arise in every numerical computation.
Here we improve the algorithm for computing the spectral gap and provide several new lower bounds for $\kappa(\SL(n,\Z ), E(n))$ and $\kappa(\SL(n,\mathbb{F}_p), E(n))$, where $E(n)$ is the set of elementary matrices (i.e. matrices which differ from identity by $\pm 1$ on a single entry off the diagonal). Related results were recently obtained independently by Fujiwara and Kabaya \cite{Fujiwara2017}.

The paper is organised as follows: in the remaining part of the introduction we provide the necessary background and notation. In Section 2 we provide a detailed description of the algorithm, Section 3 is devoted to implementation details, and in the last section we gather the results of numerical experiments.

\paragraph{Acknowledgements.}
We would like to thank Dawid Kielak and Taka Ozawa for interesting comments. 

\paragraph{Group Laplacian.}
Let $G = \langle\;\! S\mid \!R\;\!\rangle$ be a finitely presented group, with finite generating set $S$ and finite set of relations $R$.
Usually the set $S$ is assumed to be symmetric, i.e. $S^{-1} = S$, where $S^{-1}$ denotes the set of point-wise inverses of $S$.
The \emph{unnormalized Laplace operator}\footnote{Note that some authors (\cite{Ozawa2016, Bekka2008}) use normalised Laplace operator, (which differs from $\Delta_S$ by $\frac{1}{|S|}$ normalisation constant), while others (\cite{Netzer2015, Fujiwara2017}), use $\Delta_S$ as defined here.} $\Delta_S\in \mathbb{R}[G]$, associated to $S$, is defined as 
\[\Delta_S = \frac{1}{2}\sum_{g\in S} (1-g)^*(1-g) = |S| - \sum_{g\in S} g,\]
where $^*$ is the standard involution on $\mathbb{R}[G]$ given by $g^* = g^{-1}$ for $g\in G$ and $a^* = a$ for $a\in \mathbb{R}$. Throughout the paper we will drop the subsript $S$ from $\Delta$, whenever 
it does not lead to confusion.
Spectral properties of $\Delta_S$ contain a great deal of information on $(G,S)$, and thus have been intensely studied.

\paragraph*{Kazhdan property (T) and spectral gap.}
For $\pi\colon G\to \mathcal{B}(\Hilbert)$, an orthogonal representation of $G$ on a 
(real) Hilbert space $\Hilbert$ denote by $\Hilbert^\pi=\{ v\in \Hilbert\colon \pi_gv=v \text{ for all } g\in G\}$ the 
(closed) subspace of $\pi$-invariant vectors.
We define \[\kappa(G,S,\pi) = \inf\left\{\sup_{g\in S}\|\pi(g)\xi - \xi\|_\Hilbert\colon \xi\in \left(\Hilbert^\pi\right)^{\perp}, \|\xi\|=1\right\}.\]
The Kazhdan constant $\kappa(G,S)$ is defined as the infimum of $\kappa(G,S,\pi)$ over all orthogonal representations $\pi$ of $G$.
We say that $G$ \emph{has Kazhdan property (T)} if and only if there exists a finite generating set $S$ such that $\kappa(G,S) > 0$.
Proving that a group $G$ has property (T), and further estimating the Kazhdan constant of $(G,S)$ is usually a difficult task. 

Analogously to $\kappa(G,S)$ we define $\lambda(G,S,\pi)$ as the spectral gap of $\pi(\Delta_S)$; i.e., its first non-zero eigenvalue, and $\lambda(G,S)$ 
as the infimum of $\lambda$s over all representations $\pi$ without non-zero invariant vector.
The spectral gap is related to Kazhdan constant of $(G, S)$ by the inequality (see \cite[Remark 5.4.7]{Bekka2008})
\[\sqrt{\frac{2\lambda(G,S)}{|S|}} \leq \kappa(G,S).\]

The following characterization of property (T) is instrumental in our further considerations.
\begin{theorem}[Ozawa, \cite{Ozawa2016}]
A discrete group $G = \langle S\,|\,\ldots\,\rangle$ has Kazhdan's property (T) if and only if there exists constant $\lambda > 0$ and $\xi_1, \ldots, \xi_k\in \mathbb{R}[G]$ such that 
\begin{equation}
\Delta_S^2 - \lambda \Delta_S = \sum_{i=1}^{k} \xi_i^*\xi_i.\label{eq:SOS-decomposition}
\end{equation}
\end{theorem}
\noindent It is clear that $\lambda$ in the theorem is a lower bound for $\lambda(G,S)$, and in fact the maximal such $\lambda$ is equal to $\lambda(G,S)$, although we will not use the equality.

The possibility of obtaining certified (i.e. mathematically correct) estimate on $\lambda$ by numeric optimisation is a consequence of a lemma which controls residuals bound to occur in the computer-based calculations. Let $\omega[G]$ denote the kernel of the augmentation homomorphisms $\varepsilon\colon \mathbb{R}[G] \to \mathbb{R}$, and for $V<\mathbb{R}[G]$ let $V^h$ denote the subspace of $^*$-invariant elements of $V$.
\begin{lemma}[Ozawa]\label{lem:Ozawa_Delta_order_unit}
Suppose that $r\in \omega[G]^h$. Then $r + R \Delta \geq 0$, i.e. admits a sum-of-squares decomposition for $R$ large enough.
\end{lemma}
Netzer and Thom (\cite[Lemma 2.1]{Netzer2015}) gave an explicit bound on $R$ in terms of $wl_S$, the word length in $S$:
\[R\leq 2^{m}\|r\|_1, \]
where 
\[m = m(r,S) = \max_{g\in \supp(r)}2wl_S(g) - \chi(S),\]
and $\chi(S) = 1$ if $S$ contains elements of order $2$, or $\chi(S) = 2$ otherwise.
% This is a natural estimation that arises in the proof, however, for computational purposes we propose a better (if less readable) bound on $R$:
% \begin{lemma}\label{lem:residual-estimation}
%  Suppose tha $r \in \omega[G]^h$. Decompose $r$ as $\sum_g \frac{r_g}{2}(2- g^* - g)$. Then 
%  \[R \leq \sum_{g\neq e, r_g < 0} 2^{m(g)}|r_g|,\]
%  where $m(g) = 2wl(g) - \chi(\supp(g))$.
% \end{lemma}
% 
% The lemma follows directly from the proof of \cite[Lemma 2.1]{Netzer2015}, by estimating every summand of $r$ directly. We note that the numerical effectiveness of the lemma above is more in keeping $\chi(\supp(g))$ than $|r_g|$ separated.

\section{The algorithm.}

We provide only a bare minimum in the paragraphs below, as the theory and the approach has been described in \cite{Netzer2015}. For general theory of convex (conic) optimisation, we refer to excellent book \cite[Chapter 4]{Boyd2004}, or the survey \cite{Vandenberghe1996}. For the thorough discussion of the method explained below we suggest the very well written survey \cite{Netzer2016}.

\paragraph{Sum of squares and semidefinite programming}Let $\mathbf{x}$ be a fixed, ordered basis of a finite dimensional subspace $V\subset\mathbb{R}[G]$. A sum-of-squares decomposition in $V$ (as in (\ref{eq:SOS-decomposition})) is equivalent to the existence of semi-positive definite matrix (the so called Gramm matrix) $P$ satisfying
\[\Delta^2 - \lambda \Delta = \mathbf{x}^*P\mathbf{x}^T.\]
Indeed, by semi-positive difiniteness of $P$ we can find a square matrix $Q$ such that $P = QQ^T$.
Matrix $Q = [q_1|\ldots|q_k]$ can be seen as change-of-basis operator, as setting $\mathbf{x}q_i = \xi_i$ we obtain 
\[\mathbf{x}^*P\mathbf{x}^T = (\mathbf{x}Q)^*(\mathbf{x}Q)^T =\sum_{i=1}^k\xi_i^*\xi_i.\]

Semi-definite optimisation is a subclass of conic optimisation, where (in the primal form) the objective function is a linear functional minimised over a spectrahedron (i.e. the intersection of an affine space and the cone of semi-positive definite matrices).
The optimisation problem for our sum-\hspace{0pt}of-\hspace{0pt}squares decomposition can be described as a semidefinite optimisation problem (SDP)  

\begin{equation}
\begin{aligned}
  \text{minimise} & & -\lambda\\
\text{subject to} & & \mathbf{x}^*P\mathbf{x}^T &= (\Delta^2 - \lambda\Delta)\\
                  & & \lambda &\geq 0\\
                  & & P &\succeq 0,
\end{aligned}\label{eq:SDP-problem}
\end{equation}
where $P\succeq 0$ denotes semi-positive definetness.
Note that, in this formulation, entries of $P= (p_{ij})$ and $\lambda$ are optimised variables, while the equality constraint defines a set of linear constraints involving $p_{ij}$ and $\lambda$. 
Providing just one solution $(\lambda, P)$ of the problem certifies property (T) for the group $G$ (if $\lambda > 0$), however it is known that deciding feasibility (non-emptiness of the set of solutions) of an SDP is an NP-co-NP-hard problem \cite[Theorem 2.30]{Tuncel2010}. Nonetheless, there are specialised languages (\emph{algebraic modeling languages}) to specify SDP problems and softwares (\emph{solvers}) which obtain \emph{provably approximate} values $\lambda_0$ and $P_0$ with surprising efficiency.

The two key problems here are:
\begin{itemize}
\item the choice of subspace $V$, and
\item turning the approximate solution $(\lambda_0, P_0)$ into exact one, i.e. dealing with the residual of the numerical solution \[r = \Delta^2 -\lambda_0 \Delta - \mathbf{x}^*P_0 \mathbf{x}\approx 0.\]
\end{itemize}

\paragraph{The choice of $V$.} For classical sum-of-squares problem, in the (Laurant) polynomial algebra, the choice of $V$ is straightforward. Given a polynomial $f$ of (absolute) degree $2d$ in $k$ variables, it is enough to take $V$ spanned by all monomials of (absolute) degree not exceeding $d$. However, in the case of optimisation in the group algebra $\mathbb{R}[G]$ it is not clear what is the optimal $V$, as the set of relations $R$ plays an important role in this choice. 

Since $\Delta^2$ is supported on $B_2(e,S)$ it seems natural to include $\langle S\rangle$ in $V.$ Feasibility of the problem for $V = \langle B_1(e,S) \rangle$ is a decomposition as sum-of-squares of elements of length $2$ (and corresponds to the 1-st filtration in Lasserre's hierarchy of polynomial optimisation \cite{Lasserre2001}). Taking $V=\langle B_2(e,S)\rangle $ or $V=\langle B_3(e,S)\rangle$ (length $\leq 4$ or $\leq 6$) may give additional degrees of freedom in the optimisation, and thus lead to better (in terms of $\lambda$) solution to problem (\ref{eq:SDP-problem}).

Indeed, in most cases the choice of larger $V$ proved to be beneficial: see Tables~\ref{tab:computation-metrics-SL(2)} and \ref{tab:computation-metrics-SL(3-4-5)}.
However, to specify constraints on the coefficients of $P$ for $V=\langle B_d(e,S)\rangle$, we need the full multiplication table on $B_d(e,S)$, which requires generation of $B_{2d}(e,S)$.
Due to the growth of $G$, the dimension of $V$ grows exponentially with $d$, thwarting any attempts of numerical solution, except for small groups.
Moreover, the error term grows exponentially with $d$ as well obstructing bound certification, even if the numerical optimisation was successful.
In the last section we report our experiments for $V=\langle B_d(e,S)\rangle$ with $d \leq 5$, i.e. the decomposition, if obtained, is a  sum of squares of elements of length at most~$5$.

We note also a basic fact for $\SL(2,p)$: for every $d$ there exists a prime $p$ such that no decomposition 
\[\Delta^2 - \lambda\Delta = \sum \xi_i^*\xi_i\] is attainable for $\xi_i \in \langle B_d(e,S) \rangle$.
In particular, for $p > 7$ it is not possible to obtain a positive estimate on $\lambda(\SL(2,p), E(2))$ by a sum-of-squares decomposition of $\Delta^2 - \lambda\Delta$ over $\langle B_2(e,S) \rangle$.
Indeed, for a fixed $d$, the ball $B_{2d}(e,S)$ is mapped isometrically by the modular projection $\SL(2,\Z) \to \SL(2,p)$ for $p$ large enough.
Thus a decomposition of $\Delta^2 - \lambda\Delta$ obtained for $\SL(2,p)$ would lift to a similar decomposition for $\SL(2,\Z)$, which is a contradiction. This (partially) explains why a positive estimate for $\SL(2,7)$ over $\langle B_2(e,S) \rangle$ is not obtainable by numerical optimisation -- the balls of radius $2$ in $\SL(2, \Z)$ and $\SL(2, 7)$ are very similar (they differ by only $4$ elements).
On the other hand, the optimisation finishes successfully over $\langle B_3(e,S)\rangle$ and $\langle B_4(e,S)\rangle$ (which is all of $\mathbb{R}[\SL(2,7)]$).

If $G$ has property~(T) and $V = \mathbb{R}[G]$, then problem (\ref{eq:SDP-problem}) has solutions for all $0 \leq \lambda <\lambda(G,S)$ and is thus strictly feasible. By Slater's condition this implies no duality gap and therefore the uniqueness of the optimum value (see \cite[Chapter~5]{Boyd2004}). On the other hand, the spectral gap $\lambda(G,S)$ can be expressed as
\[\lambda(G,S) = \sup\left\{\lambda \colon \Delta^2 - \lambda\Delta =\sum\xi_i^* \xi_i \text{ for some }\xi_i\in \mathbb{R}[G]\right\},\]
so in the case of finite group $G$ the supremum is clearly attained on $V$.
For a few first primes $p$, we can extend $V$ to the whole of $\mathbb{R}[\SL(2,p)]$ and therefore obtain a solution which should be relatively close to the actual $\lambda(\SL(2,p), E(2))$.

\paragraph{The residual.}

Although our treatment of the residual is similar to that of \cite{Netzer2015}, we provide a short description as we will modify and explain some of the details in the next section.

In numerical optimisation we are always left with residuals, i.e. for the problem of the decomposition of $\Delta^2 - \lambda \Delta$ into sum-of-squares, a numerical solution $(\lambda_0, P_0)$ is a floating-point approximation
\[\Delta^2 - \lambda_0 \Delta - \mathbf{x}^*P_0\mathbf{x} = r, \quad \|r\|_1 \leq \varepsilon.\]
To turn the numerical solution to problem (\ref{eq:SDP-problem}) into exact one, i.e. provide a certified bound on $\lambda$, we have to modify the solution as follows.
\begin{enumerate}
\item Compute $Q = \operatorname{Re}(\sqrt{P_0})$: the real part of square root of $P$;
\item Approximate $Q$ by $Q_{\mathbb{Q}}$: project the floating-point matrix $Q$ to a close, rational-valued one;
\item Correct $Q_{\mathbb{Q}}$ to $Q^{\omega}_{\mathbb{Q}}$: each column $q_i$ of $Q_{\mathbb{Q}}$ is projected to $q_i^\omega$ which represents an element of the augmentation ideal, i.e. $\mathbf{x}q_i^\omega\in \omega[G]$;
\item Compute the decomposition $\sum \xi_i^*\xi_i = \mathbf{x^*}Q^{\omega}_{\mathbb{Q}}{\left(Q^{\omega}_{\mathbb{Q}}\right)}^T \mathbf{x}^T$ and the residual: \[r = \Delta^2 - \lambda_0 \Delta - \sum\xi_i^*\xi_i.\]
\item Obtain certified bound for the spectral gap: 
\begin{align*}
\Delta^2 - \lambda_0 \Delta = \sum \xi_i^*\xi_i + r &\geq \sum \xi_i^*\xi_i + 2^{m(S)}\|r\|_1 \Delta,
\intertext{hence}
\Delta^2 - \left(\lambda - 2^{m(S)}\|r\|_1\right)\Delta & \geq \sum \xi^*_i\xi_i \geq 0,
\end{align*}
which certifies $\lambda(G,S) \geq \lambda_0 - 2^{m(S)}\|r\|_1.$
\end{enumerate}

\section{Implementation details.}

It is clear that the difference between numerical $\lambda_0$ and the certified lower bound for $\lambda(G,S)$ grows with the approximation error.
Thus, it may be favourable not to simply maximize $\lambda$, but rather turn problem (\ref{eq:SDP-problem}) into multi-objective optimisation with minimisation of $\|r\|_1$ as the second objective.
We did not pursue this direction, as the objective function
\[P\longrightarrow \lambda - 2^{m(S)}\big\|r\big(Q_{\mathbb{Q}}^\omega\big)\big\|_1\]
outlined above is neither linear, nor differentiable, nor convex, making it very hard - if not impossible - to execute efficiently an optimisation problem.
We do, however, introduce a few changes to the problem which will exploit the strong dependence of the certified bound on $\|r\|_1$.
These changes lead to significant increase in accuracy (hence better certified bound) while at the same time (surprisingly) decrease the time required to obtain the solution.

\paragraph{Problem modification.}
Several types of solvers are available for semidefinite optimisation problems. 
Interior-point solvers are usually very accurate and converge quickly (in terms of iterations) on small to middle sized problems, but require substanital 
computing resources.
We used first-order solver (based on \emph{alternating direction method of multipliers} -- ADMM) which, while achieving moderate-to-poor accuracy \cite{He2017}, obtains an approximate solution $(\lambda_0, P_0)$ of problem (\ref{eq:SDP-problem}) relatively quickly.
Moreover, compared to interior-point solvers, first-order solvers are memory efficient, which allows us to tackle much larger problems.
At first glance, accuracy provided by an ADMM solver may seem counter-productive to our purposes.
However, after obtaining an approximate solution, we modify the original problem (\ref{eq:SDP-problem}) by constraining additionally both $\lambda$ and $P = (p_{ij})$ and run the solver again on following problem (with increased requested accuracy).
\begin{equation}
\begin{aligned}
% \begin{split}
  \text{minimise} & & -\lambda\\
\text{subject to} & & \mathbf{x}^*P\mathbf{x}^T &= (\Delta^2 - \lambda\Delta)\\
                  & & \lambda &\geq 0\\
                  & & \lambda & \leq (1-\delta)\lambda_0\\
                  & & P &\succeq 0\\
                  & & \sum_{i,j} p_{ij} & = 0.\\
% \end{split}
\end{aligned}\label{eq:SDP-problem2}
\end{equation}
% We provide the obtained solution $(\lambda_0, P_0)$ as a warm-start to the solver, and increase the requested accuracy.
Constraining $\lambda$ from above forces the solver to increase accuracy of $P$ rather than over-estimating $\lambda$ ($\delta$ is a chosen small positive number).
The additional constraint on the sum of elements in $P$ minimises the projection distance introduced in step~$3$: $q_i \mathbf{x}\in \omega[G]$ implies $\sum_{ij}p_{ij}=0$.  Note that although the constraint was implicit in the original problem, making it explicit to the solver boosts the accuracy of the (constrained) solution even further.

\paragraph{Validated Numerics.}
Netzer and Thom \cite{Netzer2015} (as well as Fujiwara and Kabaya \cite{Fujiwara2017}) performed the whole validation step in rational (i.e. exact) arithmetic. However, calculations with rational numbers are expensive both in time and space: additions require performing divisions and cause numerators and denominators grow at an exponential rate. Especially step $4$ is problematic as coordinates of $r$ get close to $0$, but the space required to store $r$ grows exponentially with the dimension of $V$. Indeed, for large $S$, these steps take longer than the actual optimisation phase. While steps $2$ and $3$ must necessarily be performed in exact rational arithmetic to apply Lemma~\ref{lem:Ozawa_Delta_order_unit} (or its quantitative version), there is no need to continue calculations with rationals afterwards.

After making suitable adjustments of $Q$ in rational arithmetic, we turn in step 4 to \emph{interval arithmetic}~\cite{Tucker2011}.
This technique of computation does not produce as result an exact number, but a real interval together with \emph{mathematical guarantee} that the exact result belongs to the interval. 
All calculations are carried out on real intervals, and these intervals are propagated through every function. 
Such computations are much easier for computer to perform than rational arithmetic, as they consist of floating-point operations and \textit{directed rounding}.
As the the aim of the procedure is to compute the lower bound for $\lambda(G,S)$, the left end of the resulting interval can be used as such a bound.
While interval arithmetic is inherently less precise than rational arithmetic, we can afford the loss of precision in the step due to high accuracy of the solver.

\section{Computed bounds}

We report on the obtained numerical bounds and computation metrics in Tables~\ref{tab:computation-metrics-SL(2)} and~\ref{tab:computation-metrics-SL(3-4-5)}.
The dominating part of the computation is either the pre-solve phase or the optimisation of the constrained SDP (\ref{eq:SDP-problem2}) (which is parallelisable), thus, in principle, bounds for larger groups could be obtained on a cluster computer.
In the pre-solve phase we stop the solver after the primal/dual objective has stabilised (usually this happened after either accuracy of $10^{-5}$, or 20000 iterations have been reached). In some cases of smaller groups it is possible to obtain high accuracy even in the pre-solve phase. The choice of $\delta$ was based on the solvers rate of convergence in the pre-solve phase in an ad-hoc manner.
The certified bound is not stable with respect to $\delta$: e.g. for problem of $\SL(3,\Z)$ with $d=2$, changing the upper bound on $\lambda$ to $0.28$ resulted in solvers stall and overall lower certified bound due to larger norm of the residual.
This may serve as an indicator that $\lambda_0$ exceeds the maximal $\lambda$ attainable on $V$.

An interesting case present numerical bounds for group $\SL(2,19)$ in Table~\ref{tab:computation-metrics-SL(2)}. 
For $d=5$ we obtain certified lower bound on the spectral gap of $0.245...$.
Even the increase of $d$ to $7$, does not yield a better numerical estimate (the certified bound becomes worse, as the SDP problem grows rapidy and the cost of the residual is $2^{4}$-times larger compared to $d=5$).
Note that products of elements in $B_{14}(e,E(2))$ saturate the whole group $\SL(2,19)$, hence no better decomposition of $\Delta^2 -\lambda\Delta$ into sum of squares can be found by semidefinite programming.
% This is surprising, as the Selberg $1/4$-conjecture would predict otherwise, that 

We tried to employ our software to $\operatorname{SOut}(\mathbb{F}_4)$ (generated by transvections) for $d=2$: this resulted in a relatively large optimisation problem with $n=1959211$ variables and $m=5937302$ constraints.
We compared the solvers behaviour to the problems of $\SL(5,\Z)$ ($n = 1256642, m = 1756344$) and $\operatorname{SOut}(\mathbb{F}_3)$ ($n = 231362, m = 131300$) which does not have property~(T). 
While the solver on $\SL(5,\Z)$ achieved convergence (in 15 hours), on $\operatorname{SOut}(\mathbb{F}_3)$ the solver stalled, with accuracy oscilating around $5\cdot 10^{-5\pm 1}$ and ever decreasing optimisation variable.
Note that this is different behaviour from the one observed for $\SL(3,\Z)$ and $\lambda_0 = 0.28$, where the solver was unable to increase the accuracy, while maintaining the value~of~$\lambda$.
The solver on problem for $\operatorname{SOut}(\mathbb{F}_4)$ exhibited a similar behaviour to $\operatorname{SOut}(\mathbb{F}_3)$ (stall, even with very small $\lambda_0$) for $\operatorname{SOut}(\mathbb{F}_3)$, i.e. we were not able to get certified positive~$\lambda$.
This suggests that $\Delta^2 - \lambda\Delta$ is not decomposable (for any $\lambda$) into sum of squares of elements in $\langle B_2(e,S)\rangle$. We consider this a (heuristic) argument supporting a common belief that $\operatorname{Aut}(\mathbb{F}_4)$ does not have property (T).

\paragraph{Software details.}
The source code used to perform computations is licensed under MIT license and accessible at \url{https://git.wmi.amu.edu.pl/kalmar/}. The computation is divided into three steps. 
\begin{itemize}
\item Computing basis of $V = \langle B_d(e,S) \rangle$ and the multiplication table;
\item Defining and solving SDP problems (\ref{eq:SDP-problem}) and (\ref{eq:SDP-problem2});
\item Estimating the certified bound on $\lambda(G,S)$.
\end{itemize}
These steps are mostly independent of each other, and it is possible to e.g. supply basis and the multiplication table computed in other programs (e.g. GAP \cite{GAP4}) and plug it into SDP routines of ours. A slight enhancement to the implementation would be to provide (the primal and the dual) solutions of the pre-solver as a warm-start for the constrained problem. In case of large problems (e.g $\SL(5,\Z)$, or $\SL(5, p)$) this may decrease the time needed to solve the constrained SDP and thus obtain a better bound in the same time-frame.

To credit the developers we provide a few details of the softwares used. The implementation is written in \textsc{Julia} \cite{Bezanson2017} programming language focused on scientific computations. The SDP solver employed is \textsc{SCS} \cite{Donoghue2016} through the \textsc{JuMP} modelling language \cite{Dunning2015}.

% \small
% \bibliographystyle{hyperplain}
\bibliographystyle{habbrv}

\bibliography{SDPPropertyT}{}

\newgeometry{bottom=4.15cm, top=4.15cm, left=2cm, right=2cm}
\begin{landscape}
% \setlength{\tabcolsep}{6.2pt}
% \newgeometry{left=3cm,bottom=0.1cm}
\begin{table}[htb]
\begin{center}

\begin{filecontents*}{table_new.csv}
Group;S;Size;d;B2d;n;m;t-pre;lambda0;(1-d)lambda0; t-SDP;iters; prec;resid;kappa-a;kappa-n;kappa-our;upper-bound
#;;;;;;;;;;;;;;1/(42*sqrt(n)+860) (SL(n,Z));(Fujiwara);;
#;;;;;;;;;;;;;;1/(31*sqrt(n)+700) (SL(n,p));;;
\SL(2, 3);4;24;2;24;92;208;0.1;1.26794;1.2679;35;257740;0.0000000000001;0.000000000000798347;0.00134437401310125;0.7961;0.796222330758409;1
%;;;;;;;;;;;;;;;;;
\SL(2, 5);4;120;2;91;154;340;54;0.13672;0.13;290;1956900;0.00000004;0.000477794;;0.258;0.259623558253098;
;;;3;120;947;2014;1.9;0.76393;0.76393;;5760;0.000000001;0.0000000161737;;0.6145;0.618033065547791;
%;;;;;;;;;;;;;;;;;
\SL(2, 7);4;336;2;111;154;419;;;;;;;;;;;
;;;3;326;1129;2584;2.38;0.58578;0.58578;;3600;0.000000003;0.0000000750102;;0.5387;0.541192570549892;
;;;4;336;6217;12770;20.9;0.58578;0.58578;;8100;0.000000002;0.0000000382189;;;0.541191995501781;
%;;;;;;;;;;;;;;;;;
\SL(2, 11);4;1320;3;532;1129;2790;;;;;;;;;;;
;;;4;1234;6671;14576;64.9;0.38196;0.38196;;11900;0.000000006;0.000000188746;;;0.437005671734361;
;;;5;1320;34192;69704;1790;0.38196;0.38196;;46400;0.000000005;0.000000211903;;;0.436981548713444;
%;;;;;;;;;;;;;;;;;
\SL(2, 13);4;2184;3;556;1129;2814;;;;;;;;;;;
;;;4;1622;6671;14964;181;0.32486;0.32486;;26900;0.00000002;0.000000141256;;;0.403020433486939;
;;;5;2184;34717;71618;4290;0.32487;0.32486;417;24500;0.000000006;0.000000338451;;;0.40297230087439;
%;;;;;;;;;;;;;;;;;
\SL(2, 17);4;4896;4;2268;6671;15610;;;;;;;;;;;
;;;5;4690;34717;74124;905;0.29072;0.29072;;52500;0.000000008;0.000000385443;;;0.381196352679298;
;;;6;4869;154847;314590;6750;0.29072;0.29072;;106900;0.000000005;0.000000144592;;;0.381163962614516;
%;;;;;;;;;;;;;;;;;
\SL(2, 19);4;6840;4;2144;6671;15487;321;0.10781;0.107;147;40960;0.000000001;0.00000000417941;;;0.231300379936739;
;;;5;5428;34717;74863;690;0.24539;0.24539;;50700;0.000000008;0.00000032573;;;0.35021893518198;
;;;6;6830;154847;316524;12300;0.24539;0.24539;7530;128300;0.000000007;0.000000547578;;;0.349878030982227;
;;;7;6840;625522;1257885;23700;0.24539;0.24539;390;1600;0.0000001;0.0000028543;;;0.341832332584266;
%;;;;;;;;;;;;;;;;;
\SL(3, 2);6;168;2;142;497;1136;0.8;1.5858;1.5857;111;276220;0.00000000003;0.00000000225473;0.00132679915701026;;0.727025900944425;0.816496580927726
;;;3;168;10154;20476;36.6;1.58578;1.58578;;11280;0.000000001;0.000000030249;;;0.727044022287509;
%;;;;;;;;;;;;;;;;;
\SL(3, 3);12;5616;2;2278;5996;14272;13;2.70849;2.708;1060;168820;0.00000000001;0.00000000329099;;0.6716;0.671813464540324;
;;;3;5610;177311;360232;581;2.70849;2.70849;;7840;0.000000001;0.000000199957;;;0.671873847247135;
%;;;;;;;;;;;;;;;;;
\SL(3, 5);12;372000;2;5119;7382;19884;91;1.50201;1.498;398;73480;0.00000000003;0.0000000121423;;0.4981;0.499666547376148;
;;;3;100398;379757;859912;5410;1.76749;1.7674;448;2280;0.0000002;0.00000354109;;;0.542731232219656;
%;;;;;;;;;;;;;;;;;
\SL(3, 7);12;5630688;2;5443;7382;20208;94;0.74233;0.7385;820;145500;0.00000000001;0.0000000047;;0.3508;0.350832339156944;
;;;3;145812;390287;926386;7830;1.47272;1.471;31500;89700;0.000000002;0.000000103905;;;0.495142796830706;
%;;;;;;;;;;;;;;;;;
%\SL(3, 11);12;212427600;2;;;;;;;;;;;;;;
\SL(3, 11);12;;3;154278;390287;934852;19700;0.75341;0.74;76200;112440;0.000000001;0.000000383843;;;0.3511870008775;
%;;;;;;;;;;;;;;;;;
\SL(3, \Z);12;;2;5455;7382;20220;420;0.28174;0.2799;2140;380240;0.000000000001;0.00000000086;0.00107210307674996;0.2155;0.21598610933692;
;;;3;154446;390287;935021;28200;0.54197;0.5405;919000;1080460;0.000000001;0.00000052478;;;0.300136525190343;
%;;;;;;;;;;;;;;;;;
\SL(4, 2);12;;2;;5996;14702;240;2.22794;2.227;430;74560;0.00000000001;0.0000000032;0.00131233595800525;;0.609234488845086;0.707106781186548
\SL(4, 3);24;;;;83846;216472;1990;4.30793;4.3;7000;61880;0.00000000003;0.000000034;;0.5974;0.598609490400461;
\SL(4, 5);;;;;93962;261826;1980;2.80353;2.8;22600;212060;0.00000000003;0.0000000311102;;0.4812;0.483045890866183;
\SL(4, 7);;;;;93962;263098;2010;1.72087;1.72;42800;424380;0.00000000003;0.0000000205032;;0.3284;0.37859388915259;
\SL(4, \Z);;;;;93962;263122;1790;1.31686;1.315;19100;220100;0.00000000003;0.000000052;0.00105932203389831;0.3285;0.331033732793956;
%;;;;;;;;;;;;;;;;;
\SL(5, 2);20;;;;39622;101776;2060;2.92148;2.921;14900;445600;0.000000000003;0.00000000525157;0.00129985241541547;;0.540462764431669;0.632455532033676
\SL(5, 3);40;;;;584822;1556596;30800;6.00595;6;182000;143380;0.000000001;0.0000231037;;;0.547722293825758;
\SL(5, 5);;;;;628882;1754106;30000;4.24441;4.24;418000;336700;0.00000002;0.000137278;;;0.460432712809375;
\SL(5, 7);;;;;628882;1757426;23500;3.04647;3.04;722000;640800;0.00000002;0.0000884447;;;0.389870354658122;
\SL(5, \Z);;;;;628882;1757466;22400;2.6877;2.65;54900;68700;0.00000002;0.0002;0.00104831159164866;;0.36400205906011;
\end{filecontents*}
\pgfplotstabletypeset[
    skip rows between index={19}{\pgfplotstablerows},
    columns={Group,Size,d, B2d,n,m,t-pre,lambda0,(1-d)lambda0,t-SDP,prec,resid,kappa-n,kappa-our}
]
{table_new.csv}
\caption{Computation metrics for $\SL(2,p)$: the size of generating set $S$ is $4$; $|B_{2d}| = |B_{2d}(e,S)|$ is the size of support of $\mathbf{x}^*P\mathbf{x}^T$;
$n$ and $m$ are the number of variables and constraints, respectively, in SDP problems;
$t_\textsl{pre}$ and $t_\textsl{SDP}$ are the time spent by solver in the pre-solve phase and in the optimisation of the constrained problem (in seconds); in case the unconstrained solver achieved decent precision in the pre-solve optimisation we set $\delta=0$ and $t_\textsl{SDP} = t_\textsl{pre}$; the times reported are for Intel Sandy-Bridge 4-core desktop processor;
\textsl{prec} is the numerical precision obtained by the solver; blank lines are left when solver did not achieve convergence;
$\|r\|_1$ is the upper bound on the $\ell_1$-norm of the residual;
$\kappa_n$ is the numerical bound obtained in \cite{Fujiwara2017};
$\kappa$ is a down-rounded numerical approximation of our certified bound -- due to numerical accuracy of the solver we certify $\lambda > (1-\delta)\lambda_0 - \textsl{prec} - 2^{m(S)}\|r\|_1$;
Previously known (analytical) lower bound: $0.0013444$ by \cite{Kassabov2005}; the upper bound for $\kappa(\SL(2,p), E(n))$ derived by Żuk (\cite{Shalom1999}) is $1$.
}
\label{tab:computation-metrics-SL(2)}
\end{center}
\end{table}
\end{landscape}
\newgeometry{bottom=3.1cm, top=3.1cm, left=2cm, right=2cm}

\begin{landscape}
\begin{table}[htb]
\begin{center}
\pgfplotstabletypeset[
    skip rows between index={0}{19},
%     skip rows between index={35}{\pgfplotstablerows},
    columns={Group,S,Size,d,B2d,n,m,t-pre,lambda0,(1-d)lambda0,t-SDP,prec,resid,kappa-n,kappa-our}
]
{table_new.csv}

\caption{Computation metrics: For description see Table \ref{tab:computation-metrics-SL(2)}. The analytical lower bounds for $\kappa(\SL(k,\bullet), E(n))$ of \cite{Kassabov2005} are: 
$0.0013268$, $0.0013123$ and $0.0012999$ (for $\SL(k,p)$) and
$0.0010721$, $0.0010593$ and $0.0010483$ (for $\SL(k,\Z)$) for $k=3,4,5$, respectively;
The upper bounds derived by Żuk (\cite{Shalom1999}) are $0.81650$, $0.70711$ and $0.63246$ for $k=3,4,5$, respectively.
}
\label{tab:computation-metrics-SL(3-4-5)}
\end{center}
\end{table}
\end{landscape}
\end{document}